\documentclass{article}

\usepackage{amssymb}
\usepackage{amsmath}
\usepackage{mathrsfs}
\usepackage{dsfont}

\usepackage[english]{babel}

\newtheorem{theorem}{Theorem}[section]
\newtheorem{lemma}[theorem]{Lemma}
\newtheorem{e-proposition}[theorem]{Proposition}
\newtheorem{corollary}[theorem]{Corollary}
\newtheorem{definition}[theorem]{Definition\rm}
\newtheorem{remark}{\it Remark\/}
\newtheorem{example}{\it Example\/}

\begin{document}

\title{On the self-$CPG$ curves and the Bj\"orling problem}
\author{Hugo Jim\'enez P\'erez and Santiago L\'opez de Medrano.}
\maketitle

\begin{abstract}
  Schwarz's solution to the Bj\"orling problem leads to an equivalence 
  class of spatial strips $S(t)=(c(t),n(t))$ which produce equivalent 
  minimal surfaces. For the particular case when the generating strip 
  $S(t)$ belongs to some plane $E$ and $c(t)$ is a symmetric curve with 
  respect to some straight line in $E$, the symmetries of the minimal 
  surface permit  us to identify another planar (geodesic) curve 
  $\tilde c(t)$ that we call the 
  \emph{CPG} curve to $c(t)$. A simple symmetric argument shows that
  \emph{self-CPG} curves produce minimal surfaces whose adjoint 
  surface contains another \emph{self-CPG} curve. We ask for minimal surfaces 
  with \emph{self-CPG} curves which are \emph{self-adjoints}.
\end{abstract}

\section{Introduction}
Schwarz's solution to the Bj\"orling problem permit us to construct a 
lot of minimal surfaces from real analytic strips $S(t)=(c(t),n(t))$, where 
$c:I\to \mathbb R^3$ is a real analytic curve and $n:\mathbb R^3\to
\mathbb R^3$ is an unitary vector field over $c(t)$ such that 
$\langle \dot c(t),n(t)\rangle \equiv 0$.
For the case when $S(t)$ is contained in some 
plane $E$, the unitary vector field $n(t)$ is recovered from the 
principal normal field $\mathfrak n(t)=\ddot c(t)/\|\ddot c(t)\|$ 
assuming that $c(t)$ is parameterized by arc lenght. 
In this situation, $c(t)$ is a plane geodesic of the minimal 
surface $X:\Omega\to\mathbb R^3$ which solves the Bj\"orling problem.

In a general context, we can consider the set of viable 
strips $\mathscr S=\{S(t)=(c(t),n(t))\}$ (see section 2.2)
and consider equivalence classes $[S(t)]$ such that for every
$\tilde S(t)\in[S(t)]$ the minimal surface $\tilde X(w)$ which
solves the Bj\"orling problem is congruent to $X(w)$. The space 
$\mathscr S$ is very big, however
we are interested in a particular class of strips, the planar strips
which posses a simple symmetry.
Suppose that the planar curve $c(t)$ has a line of 
symmetry $\mathscr L$ which intersects it perpendicularly. A simple 
analysis of the symmetries shows that $X$ will have another symmetry 
plane $E_{\mathscr L}$ which 
intersects $E$ perpendicularly along $\mathscr L$.
The plane $E_{\mathscr L}$ will contain another planar geodesic 
$\tilde c(t) \subset X$. We say that $\tilde c(t)$ is the 
\emph{conjugated perpendicular geodesic} (\emph{CPG}) to $c(t)$.
Evidently, both belongs to the same equivalence class $[S(t)]$
for $S(t)=(c(t),\mathfrak n(t))$.

In this paper we are concerned with minimal surfaces which are solutions
to the Bj\"orling problem for strips $S(t)$ whose supporting curves $c(t)$ 
are the \emph{CPG} of themselves, up to an specific rotation. 
We call them \emph{self-CPG} curves. 
We give examples of self-\emph{CPG} curves which comes from 
some classical minimal surfaces and we relate the self-\emph{CPG} 
condition with the \emph{self-adjoint} property of minimal surfaces.

\section{The Bj\"orling equivalence for planar curves}
First we recall some well-known facts from the theory of minimal surfaces.
We follow the description given by Dierkes \emph{et al.} in 
\cite{Die1}.

\subsection{Parametric minimal surfaces and geodesics}

Let $\tilde\Omega$ be an open simply connected subset of $\mathbb R^2$ and let
$X:\tilde\Omega\to\mathbb R^3$ be a mapping of class at least $C^2$ which sends 
$w=(u,v)\in\tilde\Omega$ to $X(u,v)\in\mathbb R^3$.  The image of $X$
in $\mathbb R^3$ is a minimal surface if the mapping $X$ satisfies the 
equations 
\begin{eqnarray}
  \Delta X =0 \hspace{50pt}\\
  |X_u|^2 = |X_v|^2, \qquad \langle X_u,X_v \rangle =0
  \label{eqn:min}
\end{eqnarray}
on $\tilde\Omega$, where $\Delta$ is the Laplace-Beltrami operator. In the rest of 
this document we identify the mapping with its image and we say that $X$ is a
minimal surface in $\mathbb R^3$. 

We define the \emph{adjoint surface} to $X$ on $\tilde\Omega$ 
as the surface $X^*$ which solves the \emph{Cauchy-Riemann equations}
\begin{eqnarray}
   X_u=X^*_v,\qquad X_v=-X^*_u,
\end{eqnarray}
from where we obtain that the adjoint 
surface $X^*$ to a minimal surface $X$ is also a minimal surface.
This fact permit us to state the 
problem from the complex point of view identifying 
$\mathbb C\cong\mathbb R^2$. 

Let $f:\Omega\to\mathbb C^3$ be a holomorphic mapping defined 
on the open domain $\Omega=\tilde\Omega\setminus\{Sing(f)\}$,
lets denote by $f^\prime(w)=\frac{\partial f(w)}
{\partial w}$ the derivative of $f(w)$ with respect to $w$,
and by $\langle,\rangle:\mathbb C^3\times\mathbb C^3\to\mathbb C^3$
the Hermitian inner product on $\mathbb C^3$. 
If $\langle f^\prime(w),f^\prime (w)\rangle\equiv 0$,
vanish identically on $\Omega$, the map $f(w)$ is 
called an \emph{isotropic} (complex) curve, and the real 
and imaginary components 
\begin{eqnarray}
  X(w):=\Re(f(w))\quad {\rm and }\quad X^*(w):=\Im(f(w)),
  \label{eqn:f}
\end{eqnarray}
define minimal surfaces in $\mathbb R^3$, whether or not $\Omega$ 
is simply connected.

The tangent space at any regular point $w\in\Omega$ is spanned by the 
vectors $X_u$ and $X_v$.
Additionally, at any $w\in\Omega$, the exterior product 
$X_u\wedge X_v$ does not vanish and we identify this bivector 
with its normal (perpendicular) in $\mathbb R^3$ in the traditional 
way. In a neighborhood of $w$ the \emph{unitary normal vector}
to $X$ is well defined and it is given by
\begin{eqnarray} 
   N = \frac{X_u\wedge X_v}{\|X_u\wedge X_v\|}.
\end{eqnarray}
The map $N:\Omega\to\mathbb S^2$ corresponds to the composition
$N(w):=N\circ X(w)$ and it is called the Gauss map. Since the 
image of any subset $C\subset\Omega$ in the domain of $N$ belongs
to $\mathbb S^2$ then $N(C)$ is known as the spherical image of
$X(C)$.

Two minimal surfaces $\hat X$ and $X$ are said \emph{congruents} 
 if there exist an isometry $\varphi$ and a real number $\alpha\in\mathbb R_*$
such that  $\hat X=\alpha \varphi(X)$, where $\mathbb R_*$ is the 
real multiplicative group. If $\alpha=1$, we call
 them \emph{equivalent} surfaces.

In the rest of the section the curves are parametrized 
by arc lenght. For any regular curve $c:I\to\mathbb R^3$ 
we call \emph{tangent vector} to ${\bf t}(t)=\dot c(t)$ which is a 
unitary vector, $\kappa(t)=\|\dot {\bf t}(t)\|$ is its \emph{curvature}, 
$\mathfrak n = \dot {\bf t}(t)/\kappa(t)$ its 
\emph{principal normal} and $\mathfrak b(t)=\mathfrak n(t)\times {\bf t}(t)$
its \emph{binormal}. This give us an orthonormal frame 
$\mathscr F=\left\{ {\bf t}, \mathfrak b, \mathfrak n\right\}$
over $c(t)$ from the intrinsic geometry of the curve.

Now, we consider the curve $\gamma:I\to\Omega$ such that
$c(t):= X\circ \gamma$ is parameterized by arc lenght. 
We define the \emph{normal} by ${\bf n}(t):=N(\hat c(t))$ and the 
\emph{side normal} by ${\bf s}(t):={\bf n}(t)\times{\bf t}(t)$. 
We obtain another orthonormal frame 
$\hat{\mathscr F}=\{ {\bf t}, {\bf s}, {\bf n}\}$ over $c(t)$ from the 
intrinsic geometry of $X$. 
Both frames are related by 
\begin{eqnarray*}
    \cos \theta(t)&=& \langle {\bf n}(t), \mathfrak n(t)\rangle,\\
    &=& \langle {\bf s}(t), \mathfrak b(t) \rangle.
\end{eqnarray*}
Since ${\bf t}(t)$ is an unitary vector then
${\bf n}(t)$ is a linear 
combination 
\begin{eqnarray*}
    {\bf n}(t)=\sin\theta (t)\mathfrak b(t) + \cos\theta(t)\mathfrak n(t).
\end{eqnarray*}
We define by $\kappa_g(t)=\kappa(t)\sin\theta(t)$ the \emph{geodesic
curvature} and by $\kappa_n(t)=\kappa(t)\cos\theta(t)$ the 
\emph{normal curvature} of $c(t)\subset X(w)$ for the parameter $t$.

A curve $c\subset X$ is called a \emph{geodesic} of $X$
if its geodesic curvature $\kappa_g(t)$ vanishes for all $t\in I$, it is
called an \emph{asymptotic} curve of $X$ if its normal curvature $\kappa_n(t)$
vanishes everywhere and it is called a \emph{line of curvature}
if $\dot c(t)$ is proportional to a principal direction of $X$ 
along $c(t)$, whether or not $c(t)$ is parametrized by arc lenght. 

\subsection{The Bj\"orling's problem}

Let $c:I\to\mathbb R^3$ be a real analytic curve which
admits an holomorphic extension $c(w)\subset\mathbb C^3$
and such that $\dot c(t)\neq 0$ 
almost everywhere. Over the curve $c(t)$, consider a non-vanishing 
unitary vector field $n:\mathbb R^3\to\mathbb S^2$ perpendicular to the 
tanget vector ${\bf t}(t)=\dot c(t)$, \emph{i.e.} 
$\langle {\bf t}(t),n(t)\rangle\equiv 0$.  
The couple $S(t)=(c(t),n(t))$ defines a real analytic strip 
in $\mathbb R^3$. 

Given a strip $S(t)$ as before, 
the Bj\"orling's problem concerns in to find 
a minimal surface $X:\Omega\to\mathbb R^3$ whose normal 
field $N:\Omega\to\mathbb S^2$ contains the strip $S(t)$. It means 
that $c(t)$ must belongs to $X(w)$ fullfiling the following 
properties
\begin{eqnarray}
  X(t) &=& c(t),\qquad \forall t\in I\subset\Omega,\\
  N(t) &=& n(t),\qquad \forall t\in I\subset\Omega.
  \label{eqn:strip}
\end{eqnarray}
It is immediate from conditions (\ref{eqn:strip}) and $\langle {\bf t}
(t),n(t)\rangle\equiv 0$ that $c(t)$ is a geodesic in $X(w)$.

Schwarz has proposed a solution in \cite{Sch1} (reproduced in \cite{Sch2})
using the Weierstrass representation which was generalized by the 
Cauchy-Kovalevskaya theorem. 
Schwarz's solution to Bj\"orling's problem is given by 
\begin{eqnarray}
   X(w)=\Re\left( c(w)-i\int_{w_0}^w n(z)\wedge c^\prime(z)dz\right), 
     \hspace{30pt}z,w\in\Omega\subset\mathbb C.\label{eqn:schw}
\end{eqnarray}
where $c^\prime(w)=dc(w)/dw$.

We say that $S(t)=(c(t),n(t))$ are the \emph{Bj\"orling data} for 
$X$. $\Omega$ is associated to $S(t)$ as the maximal
domain for the holomorphic extension and, in general, they 
are open domains on Riemann surfaces. 
We say that a strip $S(t)$ is \emph{viable} if there exists 
a regular parameterization of $c$ whose holomorphic extension is
defined over a punctured Riemann surface.
In particular, all the algebraic curves gives viables strips.

The space of viable strips $\mathscr S=\left\{ S(t)=(c(t),n(t))
| S(t)\ {\rm is\ viable}\right\}$ permit us to consider
 local and global parameterized curves as the same 
Bj\"orling data. Consequently the ``space'' of complete minimal surfaces
in the Euclidian space $\mathscr X=\left\{ X\subset\mathbb R^3| 
X\ {\rm is\ a\ minimal\ surface} \right\}$ will consider small open
subsets from a minimal surface and the minimal surface itself as the same 
object. We didn't have studied the implications of this consideration 
on the \emph{Schwarzian chain problem}.

We define the \emph{Bj\"orling transformation} of a strip $S(t)$ as 
the application   
\begin{eqnarray}
    \mathfrak B:\mathscr S&\to& \mathscr X\\
     S(t) &\mapsto& \Re\left( c(w)-i\int_{w_0}^{w}n(z)
      		\wedge c^\prime(z) dz \right).
    \label{eqn:map}
\end{eqnarray}
which sends the strip $S(t)$ to the minimal surface $X(w)$.

We can give a simplified strip $S(t)$ when the curve
$c(t)$ has particular properties. A classical result of O. Bonnet 
\cite{Bon1} says that it is possible to determine $X$ when the curve 
$c$ belongs to $X$ in the following cases: 
a) $c$ is a geodesic, b) $c$ is an asymptotic line, 
c) $c$ is a line of curvature,
d) $c$ is a shadow line,
e) $c$ is a perspective line.
Then consider a planar curve $c:I\to \mathbb R^3$ contained in the plane 
$E$ and the orthonormal intrinsic frame 
$\left\{ {\bf t},\mathfrak b,\mathfrak n \right\}$ over $c(t)$. Since $c(t)$
is a planar curve then the binormal vector $\mathfrak b$ coincides with
the normal $e$ to $E$ over $c(t)$. Define the normal 
$n(t)$ over $c(t)$ by 
\begin{eqnarray}
  n(t) = \mathfrak b(t) \cos \varphi(t) + 
  \mathfrak n(t) \sin\varphi(t), 
  \qquad \varphi(t)\in\left(-\frac{\pi}{2},\frac{\pi}{2}\right),
\end{eqnarray}
where ${\bf t}(t)=\dot c(t)/\|\dot c(t)\|$. 
This gives the condition
$\langle n(t), \mathfrak b(t) \rangle \equiv \cos \varphi$, for all 
$t\in I$. We obtain an analytic strip $S(t)$ whose Bj\"orling 
transformation is
\begin{eqnarray*}
    \mathfrak B(S) = \Re\left\{ c(w) -i \left( \cos \varphi(t)\mathfrak n(t)
     	+ \sin\varphi(t) \mathfrak b(t)\right)\int_{w_0}^w 
      	\|c^\prime(z)\|dz  \right\}
    \qquad z,w\in\mathbb C
\end{eqnarray*}
For $\varphi(t)\equiv\pi/2$  we obtain the classical formulation 
\begin{eqnarray}
  \mathfrak B(S) = X(w)=\Re\left( c(w) -i\mathfrak b(t)
  	\int_{w_0}^w\|c^\prime(z)\|dz \right)
   \hspace{30pt}z,w\in\mathbb C.\label{eqn:schw:plane}
\end{eqnarray}

The Bj\"orling data in expression (\ref{eqn:schw:plane}) reduces to 
$(c(t),\mathfrak n(t))$ and we write  $S(t)=(c(t))$
since the normal vector and the principal normal to the curve coincide.
When there are not way to confusion we speak about the ``Bj\"orling 
transformation of $c(t)$'' or simply ``the Bj\"orling of $c(t)$'' and
we assume that $n(t)=\mathfrak n(t)$.

\subsubsection{The Bj\"orling classes}
We say that two Bj\"orling data $S(t)$ and $\hat S(t)$ are \emph{Bj\"orling 
related} if they produce equivalent minimal surfaces. We will write 
$S\sim \hat S$ for related Borling's data. Equivalently, if the 
Bj\"orling data are given by the curves and their principal normals
then we write $c\sim \hat c$.

The uniqueness of the solution implies that we can take two
arbitrary geodesics $c, \bar c\subset X$ and 
its spherical images $n=N|_c$ and $n=N|_{\bar c}$
with regular parameterizations to produce the
Bj\"orling data $S(t)$ and $\bar S(t)$. By construction
$\mathfrak B(S)$ and $\mathfrak B(\bar S)$ are equivalent surfaces
and $S\sim \bar S$. 
In this way, we find families of  infinitelly many related Bj\"orling data.

We consider viable strips as Bj\"orling data to have a parameterization
defined in a maximal domain, which means in some punctured 
Riemann surface.
With this condition, it is an excercise to proof the following
\begin{lemma}
  $\sim$ is an equivalence relation
  \label{lem:equi}
\end{lemma}

\begin{example}
    The strips $$S(t)= \{(t,0,0), (0,\cos(t),\sin(t))\}$$ and 
    $$\hat S(t)=\{(t,0,0),(0,\cosh(t),\sinh(t))\}$$
have the helicoid as common Bj\"orling transformation, therefore 
$ S(t)\sim \hat S(t)$.
\end{example}

We can consider the classes of equivalence $[S]$ of all viable
strips $S$ such that $\mathfrak B(S)=X(w)$. We are interested in 
particular strips such that the Bj\"orling data reduce to 
planar curves.

\subsection{Schwarz's reflections and symmetries}

Schwarz discovered some interesting symmetry properties 
using expression (\ref{eqn:schw}). Such symmetries 
were used to construct a lot of minimal surfaces concatenating 
fundamental domains  of minimal surfaces whose boundary
is a composition of straight lines and/or  plane geodesics.
In order to glue two fundamental domains they must lie in the 
interior of a regular frame called a \emph{Schwarzian chain}
$\mathfrak C$. We use those symmetries for analyse the Bj\"orling
transformation of symmetric supporting curves.

A symmetry $A$ of a parametric minimal surface $X$ induce an isometry
$\alpha:\Omega\to\Omega$ such that $N\circ\alpha = \pm A\circ N$ where
$A$ is a rigid mouvement in $\mathbb R^3$. Since the spherical image
of $X$ is invariant under translations, we are interested only in
matrices $A\in O(3)$. 

Let $\tau,\lambda:\Omega\to\Omega$
be functions given by
\begin{eqnarray*}
    \tau(w)&=& \bar w\\
    \lambda(w)&=& iw,\qquad i=\sqrt{-1}
\end{eqnarray*}
and matrices $T,\Lambda\in O(3)$ given by 
\begin{eqnarray*}
    T = \left( 
    \begin{array}{ccc}
	1 & 0 & 0\\
	0 & 1 & 0\\
	0 & 0 & -1
    \end{array}
    \right), \qquad
    \Lambda = \left( 
    \begin{array}{ccc}
	-1 & 0 & 0\\
	0 & 0 & -1\\
	0 & 1 & 0
    \end{array}
    \right),
\end{eqnarray*}
which span two representations of the diedral group $D_4$ in $\mathbb C^*$ 
and $GL_3(\mathbb R)$ respectively.
We have the identities
\begin{eqnarray*}
    \tau^2=\lambda^4={\rm Id},\quad \lambda^{-1}=\tau\lambda\tau,\qquad
    T^2=\Lambda^4={\rm Id}_4,\quad \Lambda^{-1} =T\Lambda T.
\end{eqnarray*}
and in particular, $\tau$ is anticonformal and $\lambda$ is conformal.

Considering the opposite orientation of the normal field in 
the solution of Bj\"orling's problem, Schwarz obtained the same minimal
surface with the reflected domain $\bar\Omega=\left\{ \bar w | w\in\Omega
\right\}$. It has become his celebrated \emph{reflection principle}.

\begin{lemma}\label{lem:22}
  Let $X:\Omega\to\mathbb R^3$ be a nonconstant 
  minimal surface whose domain of definition $\Omega$ contains some interval
  $I$ that lies on the real axis.
  \begin{enumerate}
    \item[$\imath$)] If the curve $c(u) = \{X(u):u\in I\}$ is contained
      in some plane $E$, and if the surface $X$ intersects $E$ 
      orthogonally at $c(u)$, then $E$ is a plane of symmetry for $X$. 
  \item[$\imath\imath$)] If the image of $l(u) = \{ X(u):u\in I\}$ is contained 
      in some line $\mathscr L$, then $\mathscr L$ is a line
      of symmetry of $X$.
  \end{enumerate}
  \label{lem:sym}
\end{lemma}
We assume that the line $\mathscr L\subset \tilde X$ belongs to the z-axis and
the plane $E$ is the xy-plane.
Then $i)$ corresponds to $X\circ\tau = T\circ X$
and $ii)$ gives $X\circ\tau =-T\circ X$.

We have selected $\mathscr L\subset$ z-axis by convenience,
in order that the spherical images of $c(u)$ and $l(u)$ concide 
in $\mathbb S^2$.
In fact, they are projections of the same real curve $h:I\to\mathbb C^3$ 
with $h(u)=c(u)+il(u)\subset f(w)$.
In this case $f:\Omega\to\mathbb C^3$ is the isotropic curve 
$f(w)=X(w)+iX^*(w)$.
These relationships are contained in the next
\begin{e-proposition}\label{prop:23}
  Let $X:\Omega\to\mathbb R^3$ be a nonconstant minimal surface and assume 
  that $X^*:\Omega\to\mathbb R^3$ is an adjoint minimal surface of $X$. 
  Choose a smooth curve $\gamma:I\to\Omega$ with $\dot \gamma(t)\neq 0$ 
  except for isolated points $t_i$ in the interval $I$, and consider the curves 
  $c(t)=X\circ \gamma(t)$ and $c^*(t)=X^*\circ\gamma(t)$. 
  The following properties holds:
  \begin{enumerate}
    \item[(i)] If $c$ is a straight arc, 
      then it is both a 
      geodesic and an asymptotic line of $X$, and $c^*$ is a planar 
      geodesic of $X^*$. The curve $c^*$ lies in some plane $E$ 
      and $X^*$ intersects $E$ orthogonally along $c^*$.
    \item[(ii)] If $c$ is a planar geodesic on $X$, 
      then $c^*$ is a straight arc (and 
      hence a geodesic asymptotic line) on $X^*$. 
  \end{enumerate}
\end{e-proposition}

Assume that $c(t)\subset X(w)$ is a geodesic contained in the $XY$-plane,
then we have 
\begin{eqnarray}
 (X+X^*)(\tau w) &=& T\circ (X - X^*)(w),\qquad w\in\Omega.
\end{eqnarray}
In other words $f(\tau w)=T\circ\overline{f(w)}$ where $T$ acts 
on $\mathbb C^3$ by the diagonal action.
This result comes from the holomorphic properties of $f$.
The reader can see \cite{Die1} for the proof of Lemma 
\ref{lem:22} and Proposition \ref{prop:23}.

\begin{definition}
Suppose that $c:I\to \mathbb R^2$ has a symmetry line $\mathscr L=\mathscr 
  L(t)$
  parameterized by $\mathscr L(t)=at+b$ with 
  $a,b\in\mathbb R^3$ and $a\neq 0$. 
  We say that $c$ 
is a \emph{perpendicular symmetric curve} with respect to $\mathscr L$ if
there exist $t_0\in I$ such that $c(t_0)\in\mathscr L$ and 
$\langle \dot c(t_0), a\rangle=0$, 
We call the point $p=c(t_0)$ a \emph{symmetry vertex} of $c$.

We say that a perpendicular symmetric curve is \emph{non-degenerated}
if its normal vector $n=\ddot c/\|\ddot c\|$ does not vanishes at its
symmetry vertex.
\end{definition}

In this paper we are concerned with perpendicular non-degenerated
symmetric curves. Non-degeneracy avoids umbilical points in the 
minimal surface at the symmetry vertex of $c(t)$. The reason is that
umbilical points in minimal surfaces implies the vanishing of the
principal curvatures $\kappa_1$ and $\kappa_2$ which are necessary
in order to get perpendicular straight arcs. It is a consequence that
at umbilical points a minimal surface is not conformal to its spherical 
image. Some examples of this failure  are the 
high order element of the Enneper Family \cite{Die1} or the high genus 
Costa surfaces \cite{Cos1}.

\begin{lemma}
    Suppose that $c(t)$ is a perpendicular symmetric curve belonging to the 
    $XY$-plane. Then 
    \begin{eqnarray}
	\Lambda^2 T\circ X(w) = X(-\bar w) 
    \end{eqnarray}
    \label{lem:25}
\end{lemma}
\emph{Proof.} This is immediate from the fact that $\bar w =\tau w$
and $-w=\lambda^2 w$ then $X(-\bar w) = X(\lambda^2\tau w)$, and using Lemma
\ref{lem:22} we obtain $X(\lambda^2\tau w) = \Lambda^2T\circ X(w)$.

$\hfill\square$


\begin{lemma}\label{lem:24}
   Let $c:I\to\mathbb R^3$ be a (non-degenerated) perpendicular 
   symmetric curve and $X(w)=\mathfrak B(c)$ its Bj\"orling 
   transformation. Then $\hat c(t)=X(\lambda t)$, $t\in I$
   is a (non-degenerate) perpendicular symmetric curve.
\end{lemma}

\emph{Proof.} We suppose $c(t)\subset XY$-plane. Define $\hat c(t)
=X(\lambda t)$ which is a well defined space curve. We must prove that
$\hat c$ is a non-degenerated (planar) perpendicular symmetric curve. 
Using Lemma \ref{lem:25}
we verify that $y(-\bar w)=-y(w)$. Then $y(it)=y(\lambda t)\equiv 0$ for all 
$t\in\mathbb R$. Writting $\hat x(t)=x(\lambda t)$ and $\hat z(t)=z(\lambda t)$
we obtain 
\begin{eqnarray}
    \hat c(t) = \left( \hat x(t), 0, \hat z(t) \right).
\end{eqnarray}
Which implies that $\hat c(t)$ is a planar curve. Applying 
$X(\tau w) = T\circ X(w)$ with $w=\lambda t$
we have 
\begin{eqnarray*}
   X(\tau\lambda t) = \left( x(\lambda t), 0, -z(\lambda t) \right) 
\end{eqnarray*}
it means 
\begin{eqnarray*}
   \hat c(- t) = \left( \hat x(t), 0, -\hat z(t) \right). 
\end{eqnarray*}
Then $\hat c(t)$ is symmetric with respect to the $X$-axis.
Finally, its principal normal at the symmetry vertex does not vanish
since $\hat{\mathfrak n}(0) = -\mathfrak n(0)$ and $c(t)$ is 
non-degenerated.

We conclude that $\hat c(t)$ is a non-degenerated perpendicular 
symmetric curve and $\hat c(t)\in [c(t)]$ by construction.
$\hfill\square$

\begin{definition}
   Two perpendicular symmetric planar curves $c$ and $\hat c$ are called 
   \emph{conjugated perpendicular geodesics under the Bj\"orling 
   transformation} (or simply $CPG$), if for any parameterization of 
   $c(t)$ such that $c(t)=X(t)$, for all $t\in I$ then 
   $\hat c(t)=X(\lambda t)$ up to sign.
\end{definition}

In what follows we write only $CPG$ to mean ``the conjugated perpendicular
geodesic curves under the Bj\"orling transformation''. 

%
We recall if $c(t)$ is an algebraic 
curve its analytic version $c(z)$ will be defined in some punctured Riemann surface
and we can obtain global $CPG$ curves.

Examples of 
$CPG$ curves are the following:
\begin{itemize}
 \item The circle and the catenary: both generate the Catenoid.
 \item The parabola and the cycloid: both generate the Catalan surface.
 \item The ellipse and a class of elliptical roulette: both generate the 
     Elliptic catenoid studied in \cite{Jim3}.
 \item The cubic $(t^2, t^3/3-t)$ with itself: generate the 
     Enneper surface. 
\end{itemize}
The last example has the property that if $c(t)\subset XY$-plane
then $\hat c(t)=\Lambda\circ c(t)$, $t\in I$ 
as defined above. We call them \emph{self-CPG} curves.
In fact, if 
$c:I\to\mathbb R^3$ is a self-$CPG$ curve in the ${XY}$-plane, symmetric 
with respect to the $X$-axis and $X(w)=\mathfrak B(c)$ then 
$X (\lambda w) = 
\Lambda\cdot X(w)$.

In general, we consider the condition $X( \lambda t) = \Lambda\circ X(t)$
for $t\in I$ as the definition of the self-$CPG$ curves.

\begin{remark}
    The $CPG$ condition is not an equivalence relation. In \cite{Jim3}
    the author shows that the ellipse has two different $CPG$s, $c_1(t)$
    and $c_2(t)$, which corresponds to the vertices of the ellipse but
    $c_1$ and $c_2$ are not $CPG$ curves. The $CPG$ condition is not 
    transitive.
\end{remark}

\begin{e-proposition}\label{prop:25}
    Let $c,\hat c:I\to\mathbb R^3$ be two CPG (planar) curves and
    $X(w)=\mathfrak B(c)$ such that $c(t)$ is contained in the 
    $XY$-plane and $\hat c(t)=X(it)$ contained in the $XZ$-plane.
  Then $\hat c(t)=\Lambda\cdot c(t)$ if and only if $X(t+it)$ and $X(t-it)$ are 
  perpendicular straight lines in $X(w)$. 
\end{e-proposition}
{\it Proof.}
We begin with the necessity. We suppose $c,\hat c$ are $CPG$ and 
$X(t+it)$ and $X(t-it)$ are perpendicular straight arcs. Since
$X(0)$ is not umbilical then any neigborhood of $X(0)$ is conformal
to the disc $|z|<r$ for $z\in\Omega$ and $r>0$ small. Since $c(t)$
is contained in the $XY$-plane and $\hat c(t)$ in the $XZ$-plane, then 
$X(t+it)$ belongs to $(0,y,y)$ and $X(t-it)$ belongs to $(0,y,-y)$.

Since $X(t+it)$ is a symmetry line every point in $c=(x, y,0)$
is mapped under the symmetry to $\hat c=(-x,0,y)$. It means that
\begin{eqnarray}
    \hat c = \Lambda T\cdot c.
    \label{eqn:sym1}
\end{eqnarray}
The symmetry with respect to $X(t-it)$  implies that $c=(x,y,0)$
is mapped to $\hat c=(-x,0,-y)$. It means 
\begin{eqnarray}
    \hat c = T\Lambda \cdot c.
    \label{eqn:sym2}
\end{eqnarray}
Both curves are invariant under $T$ therefore (\ref{eqn:sym1}) and
(\ref{eqn:sym2}) implies $\hat c=\Lambda\cdot c$. Finally, 
$X:\Omega\to\mathbb R^3$ is conformal and an isometry then the
holomorphic extension preserves distances from $c(t)$ to $\hat c(t)$,
we conclude $c(t)$ is self-$CPG$.

Now the converse.
We write $t^\prime=(1-i)t$ and we have that $\lambda t^\prime=\tau 
t^\prime$. Since $c(t)$ is self-CPG we have 
\begin{eqnarray*}
    \Lambda X(t^\prime)=X(\lambda t^\prime)=X(\tau t^\prime)=TX(t^\prime),
\end{eqnarray*}
then $x(t-it)=-x(t-it)$ for all $t\in I$ and consequently $x(t-it)\equiv 0$.
Additionally we obtain $y(t-it)=-z(t-it)$ for all $t\in I$ then $X(t-it)$
is contained in the line $(0,y,-y)\subset \mathbb R^3$.

On the other hand we write $t^{\prime\prime}=(1+i)t$ and we consider the
identity $\lambda\tau\lambda\tau={\rm Id}$ to obtain 
$\lambda t^{\prime\prime} =  \lambda^2\tau\lambda\tau t^{\prime\prime} 
= \lambda^2\tau t^{\prime\prime}$. The last equality comes from the 
invariance  $(1+i)t = i\cdot \overline{(1+i)}t$. Then 
\begin{eqnarray*}
    \Lambda X(t^{\prime\prime}) = 
    X(\lambda t^{\prime\prime}) =
    X(\lambda^2\tau t^{\prime\prime}) =
    \Lambda^2 T X(t^{\prime\prime}), 
\end{eqnarray*}
equivalently $X(t^{\prime\prime})= \Lambda TX(t^{\prime\prime})$.
We obtain $x(t+it)\equiv 0$ and $y(t+it)=z(t+it)$, therefore $X(t+it)$ is
contained in the line $(0,y,y)\in\mathbb R^3$. Perpendicularity is 
obvious.
$\hfill\square$

\begin{theorem}\label{teo:main}
    Let $X:\Omega\to\mathbb R^3$ be a minimal surface such that 
    $X(w)=\mathfrak B(c)$ for a self-$CPG$ curve $c:I\to\mathbb R^3$, 
    $I\subset \Omega$.
    Then the adjoint surface $X^*:\Omega\to\mathbb R^3$ is generated 
    by another self-$CPG$ curve $c^*:I^\prime\to\mathbb R^3$ with 
    $X^*(w)=\mathfrak B(c^*)$, for $I^\prime\subset\Omega$.
\end{theorem}
{\it Proof.}
Since $c(t)$ is self-$CPG$ then it belongs to some plane
$E\subset \mathbb R^3$ which, as before, we assume $E=XY$-plane
and it is symmetric with respect to the $X$-axis. 
From Proposition \ref{prop:25},  
$X(w)$ contains two perpendicular straight arcs $X(t+it)$ and
$X(t-it)$ contained in the lines $(0,y,y)$ and $(0,y,-y)$ respectively. 

Applying
Proposition \ref{prop:23}, the $CPG$ curves $c(t)$ and $\hat c(t)$ are mapped
to two perpendicular straight arcs $X^*(t)$
and $X^*(it)$ on the adjoint surface. Meanwhile the straight arcs
$X(t+it)$ and $X(t-it)$
will be mapped to two planar geodesics $c^*(t):=X^*(t+it)$ and 
$\hat c^*(t):=X^*(t-it)$. By the self-$CPG$ definition, $c(t)$ is 
non-degenerated at $t=0$, then $\mathfrak n(0)=-\hat{\mathfrak n}(0)\neq 0$ 
and the point $X(0)$ is not umbilical. These assures that the geodesics 
$c^*(t)$ and $\hat c^*(t)$ are perpendicular and therefore they are $CPG$.

From Proposition \ref{prop:25} we have
$\hat c^*(t)=\Lambda\cdot c^*(t)$ or $\Lambda \cdot \hat c^*(t)= c^*(t)$ 
and consequently $c^*(t)$ is self-$CPG$.
%
%
%
%
%
$\hfill\square$


\begin{corollary}
    If $c:I\to\mathbb R^3$ is self-$CPG$ then the isotropic curve 
    $f(w)=X(w)+iX^*(w)$ where $X(w)=\mathfrak B(c)$ has a $D_{4}$
    symmetry.
\end{corollary}
{\it Proof.} It is enough to define $\Lambda f(w) = f(\lambda w)$ and $T\overline{ f(w)}
=f(\tau w)$ by diagonal action.
$\hfill\square$


\begin{definition}
   Let $X:\Omega\to \mathbb R^3$ be a minimal surface 
   and $X^*:\Omega \to \mathbb R^3$ an adjoint surface to $X$.   
   We say that $X$ is \emph{self-adjoint} if there exists 
   an orthogonal matrix $R\in O(3)$ and an isometry $\rho:\Omega\to\Omega$
   such that 
   \begin{eqnarray*}
   	R\circ X^*(w)=X\circ \rho (w), \qquad w\in\Omega.
   \end{eqnarray*}
\end{definition}

\begin{corollary}
  Every isotropic curve $f:\Omega\to\mathbb C^3$ whose components 
  are self-adjoint minimal surface $R X^*(w)=X(\rho w)$ obtained by the Bj\"orling 
  transformation of a self-$CPG$ curve $c:I\to \mathbb R^3$ has a $D_8$ symmetry.
\end{corollary}
{\it Proof.} This is immediate from the fact that the minimal surface $X(w)=\mathfrak B(c)$
for a self-$CPG$ curve $c(t)$ has a 
$D_4$ symmetry. Applying Proposition \ref{prop:25}, $X(w)$ posses two straight lines
which are mapped to two geodesics in its adjoint surface $X^*(w)$.
We define the complex matrix
\begin{eqnarray}
  R=
  \left( 
  \begin{array}[h]{ccc}
    i & 0 & 0\\
    0 & \frac{1}{\sqrt{2}} & -\frac{1}{\sqrt{2}}\\ 
    0 & \frac{1}{\sqrt{2}} & \frac{1}{\sqrt{2}}
  \end{array}
  \right)
  \label{eqn:Rot}
\end{eqnarray}
which maps the self-$CPG$ curves to the perpendicular straight lines
and visceversa. Writing $\rho=\exp(\pi i/4)$ then we have 
\begin{eqnarray*}
  R f(w) = f(\rho w), \qquad T\overline{f(w)}=f(\tau w)
\end{eqnarray*}
which are the generators of the $D_8$ representation in $GL(3,\mathbb C)$. Note that 
$\Lambda=R^2$ and $\lambda=\rho^2$.
$\hfill\square$

It is well-known that the Enneper surface and its adjoint are the same
geometric object. In this way, it is a self-adjoint surface. 


Another interesting example is the Costa surface 
$\mathscr C_1:S_1\to\mathbb R^3$, \cite{Cos1} where $S_1$ is the 
punctured Riemann  surface associated to $w^2=z(z^2-1)$ 
although the authors do not know a suitable parameterization 
of its self-$CPG$ supporting curve.

\section{Additional discussion}

In this note we have characterized the adjoints to minimal 
surfaces which contains self-$CPG$ curves. 
We can say that a minimal surface $X:\Omega\to\mathbb R^3$ which 
arise as the Bj\"orling transformation of a self-$CPG$ curve 
$c:I\to\mathbb R^3$, has a $D_4$ 
symmetry. This symmetry is extended to the isotropic curve 
$f(\Omega)\subset\mathbb C^3$
since the adjoint surface $X^*$ arises also as the Bj\"orling of 
another self-$CPG$ curve. 
For the case of self-adjoints surfaces coming from self-$CPG$ curves
the $D_4$ symmetry is extended to a $D_8$ symmetry on the isotropic 
curve $f(w)=X(w)+iX^*(w)$ with generator (\ref{eqn:Rot}).

%

There are other interesting examples which do not fall in the 
characterization given in this document.
 The family of algebraic curves 
   \begin{eqnarray}
     c_k(t)=\left\{\left(\frac{2}{m}t^m, 
     	\frac{1}{2m -1}t^{2m-1}-t\right)| m=4k-2,\ k\in\mathbb N\right\}, 
	\label{eqn:fam}
   \end{eqnarray}
is a family of perpendicular symmetric
curves whose Bj\"orling transformation $\mathfrak B(c_k)$ has a $D_{2k+2}$ 
symmetry. The generator 
  $\Lambda=\Lambda_k$ has the form  
  \begin{eqnarray*}
	\Lambda = \left( 
	\begin{array}[h]{ccc}
	  -1 & 0 & 0\\
	  0 & \cos \frac{\pi}{2n} & -\sin \frac{\pi}{2n}\\
	  0 & \sin \frac{\pi}{2n} & \cos\frac{\pi}{2n} 
	\end{array}
	\right), \qquad n=2k+2.
  \end{eqnarray*}
  If $k>1$ we say that $c(t)$ is a \emph{weak} $CPG$ curve. In that case
  the sufficiency condition in Proposition \ref{prop:25} is not fulfilled.

  The family (\ref{eqn:fam}) corresponds to the high order Enneper 
  surfaces \cite{Die1} and the same symmetries are shared by the high 
  genus Costa surfaces \cite{Cos1}. In both cases the symmetry arises 
  since the origin is an umbilical point.

We have several questions we are interested in to answer.

{\bf Question:} What are the conditions for some curve $c:I\to\mathbb R^2$ 
to be self-$CPG$?

{\bf Question:} Is it possible to deform the straight line to have a 
one parameter family of self-$CPG$ curves?

{\bf Question:}
There are other self-adjoint surfaces which arises from self-$CPG$ curves 
as in the case of the Enneper surface? 

{\bf Question:} 
Since there are a lot of examples
of self-$CPG$ curves whose Bj\"orling transformation gives embedded
surfaces in $\mathbb R^3$,
there exists an embedded self-adjoint surface in
$\mathbb R^3$?

\end{document}